\documentclass[orivec,envcountsect,envcountsame]{llncs}

\usepackage[a4paper, margin=1.1in]{geometry}

\usepackage[noend]{algpseudocode}
\usepackage{amsmath,amssymb}
\usepackage{algorithm}
\usepackage[normalem]{ulem}
\usepackage{enumitem}

\usepackage{amsthm}
\usepackage{graphicx}

\newtheorem{observation}[theorem]{Observation}
\numberwithin{equation}{section}

\newcommand{\cC}{\mathcal{C}}

\newcommand\R{\mathbb{R}}

\newcommand\Wout{\bgroup\markoverwith
{\textcolor{magenta}{\rule[.5ex]{2pt}{0.4pt}}}\ULon}
\newcommand\Tout{\bgroup\markoverwith
{\textcolor{blue}{\rule[.5ex]{2pt}{0.4pt}}}\ULon}
\newcommand\Mout{\bgroup\markoverwith
{\textcolor{orange}{\rule[.5ex]{2pt}{0.4pt}}}\ULon}

\newcommand{\ceil}[1]{\left\lceil #1 \right\rceil}
\newcommand{\floor}[1]{\left\lfloor #1 \right\rfloor}

\DeclareMathOperator{\boxicity}{box}
\DeclareMathOperator{\coboxicity}{co-box}

\newcommand{\Gnet}{G_{\mathrm{net}}}

\title{A polynomial-time algorithm to compute the co-boxicity of block-structured graphs}
\titlerunning{Co-boxicity of block-structured graphs}

\author{Marco Caoduro\inst{1} \and Will Evans\inst{2} \and Tao Gaede\inst{3}}
\authorrunning{Caoduro, Evans, and Gaede}

\institute{
Laboratoire d'Informatique de Paris Nord, Université Sorbonne Paris Nord, France \\
\email{caoduro@lipn.univ-paris13.fr}
\and
Department of Computer Science, University of British Columbia, Canada \\
\email{will@cs.ubc.ca}
\and
Department of Mathematics and Statistics, University of Victoria, Canada \\
\email{taogaede@uvic.ca}
}

\begin{document}

\maketitle

\begin{abstract}
    The boxicity of a graph $G$ is the minimum dimension $d$ that admits a representation of $G$ as the intersection graph of a family of axis-parallel boxes in $\R^d$.
    Boxicity has applications in ecology and operations research, where it models competition for shared resources.
    Computing boxicity is an \textbf{NP}-hard problem.
    Despite considerable theoretical interest, few graph classes are known for which the boxicity can be computed in polynomial time.
    Following the approach of Cozzens and Roberts~(1983), we study directly the boxicity of graph complements, which we call \emph{co-boxicity}.
    In this setting, the geometric intersection problem in $\R^d$ becomes a graph cover problem involving only complements of interval graphs.
    We give a polynomial-time algorithm that computes the co-boxicity of $G$ from two pieces of information about each block: its co-boxicity, and how the co-boxicity changes when certain subsets of its vertices become cut-vertices.
    Using this algorithm, we obtain polynomial-time algorithms for the co-boxicity of trees, graphs with bounded block order, block graphs, and cactus graphs.
    More generally, our method implies the union of these results: for any fixed $k \geq 3$, there is a polynomial-time algorithm to compute the co-boxicity of $G$ whenever every block of $G$ is complete, a cycle, or has at most $k$ vertices.

\end{abstract}

\section{Introduction}
\label{section-intro}
%
%
The \emph{boxicity} of a graph $G$, denoted by $\boxicity(G)$, is defined to be the minimum dimension $d$ such that $G$ can be represented as the intersection graph of a family of axis-parallel boxes in $\mathbb{R}^d$.
Boxicity was introduced in 1969 by Roberts~\cite{Rob}.
Its study was initially motivated by applications in ecology, specifically in modelling niche overlap and competition between species~\cite{1978_Cohen,1976_Roberts}.
Further applications appear in operations research, where boxicity has been used in fleet maintenance and task assignment problems~\cite{1983_Cozzens,1981_Opsut}.
\begin{example}\label{eg:boxicity-example}
    Figure~\ref{fig:boxicity-example} shows a family $\mathcal{F}$ of axis-parallel boxes in $\R^2$ with its intersection graph $G$.
    As the boxes are subsets of $\R^2$, $\boxicity(G) \leq 2$.
    Moreover, $G$ contains an induced cycle on four vertices, which cannot be represented as the intersection graph of intervals in $\mathbb{R}$.
    So, $\boxicity(G) = 2$.
    \begin{figure}
        \centering
        \includegraphics[scale=0.85]{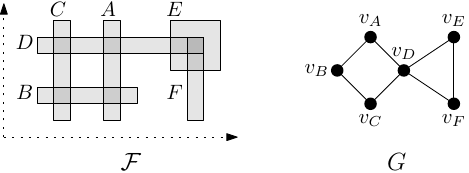}
        \caption{Family $\mathcal{F}$ of axis-parallel boxes in $\R^2$ with its intersection graph $G$.}
        \label{fig:boxicity-example}
    \end{figure}
\end{example}
Computing boxicity of an arbitrary graph is an \textbf{NP}-hard optimization problem~\cite{cozzens-phd-thesis-1981}.
While graphs of boxicity one (i.e. interval graphs) can be recognized in linear time~\cite{1975_Kellogg}, even the decision problem of whether a graph has boxicity at most $2$ is \textbf{NP}-complete~\cite{1994_Kratochvil_NP}.
Consequently, significant effort has been directed towards relating boxicity to other graph parameters (e.g.\ \cite{ABC,CS,Esp,EJ}) or computing it for specific graph classes (e.g.\ \cite{2023_Caoduro,2011_Chandran,1984_Scheinerman,1986_Thomassen}).
%

%
A \emph{co-interval graph}is the complement of an interval graph.
In 1983, Cozzens and Roberts~\cite{1983_Cozzens} initiated the study of boxicity from the perspective of covering the edges of a graph's complement using co-interval graphs.
A graph cover using only co-interval graphs is called a \emph{co-interval cover}.
Using this perspective, Cozzens and Halsey~\cite{1991_COZZENS_Coboxicity} proved that the decision problem of whether $\boxicity(G) \leq 2$ is in \textbf{P} when $G$ is the complement of a comparability graph.
Motivated by the potential of this covering perspective, we investigate the boxicity of graph complements directly.
We refer to the boxicity of the complement of a graph $G$ as the \emph{co-boxicity} of $G$, denoted by $\coboxicity(G)$.

As $\coboxicity(G) = \boxicity(\overline{G})$, co-boxicity inherits both the hardness of boxicity and the theoretical interest surrounding it.
Despite this interest, the list of known graph classes for which the co-boxicity can be computed efficiently is small.
We survey the classes we know of presently.
A first example is the class of graphs with co-boxicity at most $2$, such as complements of outerplanar graphs~\cite{1984_Scheinerman}.
In this case, computing co-boxicity is equivalent to the linear-time decision problem of whether the graph is a co-interval graph~\cite{1975_Kellogg}.
Another simple example is that of graphs of bounded order:
since every graph on $n$ vertices has co-boxicity at most $\lfloor n/2 \rfloor$~\cite{Rob}, co-boxicity can be computed in constant time by exhaustive search.
Parameterized complexity refines this by replacing the order with smaller parameters.
A result of Adiga, Chitnis, and Saurabh~\cite{2010_Adiga} shows that $\coboxicity(G)$ can be computed in \emph{fixed-parameter tractable} time when parameterized by $|V(G)| - \omega(G)$, the minimum number of vertices whose deletion leaves a clique (i.e.\ the vertex cover number of the complement).
Bruhn, Chopin, Joos, and Schaudt~\cite{2016_Bruhn} later generalized this tractability result to the minimum number of vertices whose deletion leaves $G$ complete multipartite (i.e.\ the cluster vertex deletion number of the complement).
Deciding whether $\coboxicity(G) \le 2$ is already \textbf{NP}-complete in general~\cite{1994_Kratochvil_NP}.
However, for line graphs, Caoduro and Seb\H{o}~\cite{2023_Caoduro_GD} showed that deciding whether $\coboxicity(G) \le d$  becomes tractable when parameterized by $d$.
For some classes, formulae are known so that no algorithm is required.
The path on $n$ vertices has co-boxicity $\ceil{\frac{n-1}{3}}$~\cite{1983_Cozzens}.
Similarly, Caoduro and Seb\H{o}~\cite{2023_Caoduro_GD} showed that the co-boxicity of the line graph of the complete graph on $n$ vertices is $n-2$.

Our work expands this list.
Using the co-interval cover perspective, we exploit the block structure of graphs to constrain how pieces of a co-interval cover come together.
This allows us to establish a new method that efficiently computes the co-boxicity of a graph from the co-boxicity of its blocks.

A connected graph with at least one edge is \emph{2-connected} if it has no \emph{cut-vertex}, i.e.\ a vertex whose removal disconnects the graph.
A \emph{block} of a graph is a maximal $2$-connected subgraph.
Our method reduces computing the co-boxicity of $G$ to two pieces of information about each block: its co-boxicity, and how co-boxicity changes when certain subsets of its vertices become cut-vertices.
We demonstrate the utility of our method separately on four graph classes: trees, as well as graphs with bounded block order, block graphs, and cactus graphs, i.e., graphs in which every block has at most $k$ vertices, is complete, or is a cycle or an edge, respectively.
Moreover, combining these results shows that, for every fixed $k$, the co-boxicity of a graph $G$ can be computed in time polynomial in $|E(G)|$ whenever every block of $G$ is complete, a cycle, or has at most $k$ vertices.

Before stating our theorem, we define the necessary terminology. 
A vertex is a \emph{leaf} if it has degree~$1$.
A leaf neighbour of a vertex $v$ is a \emph{pendant} of $v$.
\begin{definition}\label{def:pendant-augmentation}
  Let $Q$ be a graph and let $x : V(Q) \to \mathbb{Z}_{\geq 0}$.
  The \textbf{pendant augmentation} of $Q$ by $x$, denoted $Q^x$, is the graph obtained from $Q$ by attaching $x_v$ new pendant vertices to each vertex $v \in V(Q)$.
\end{definition}
\begin{definition}\label{def:tractable}
  A family $\mathcal{N}$ of $2$-connected graphs, each with at least three vertices, is \textbf{tractable} if it satisfies the following properties:
  \begin{enumerate}
      \item[\normalfont{\textbf{(T1)}}]\label{it:T1} for every $Q \in \mathcal{N}$ and every $u \in V(Q)$, every block of $Q - u$ belongs to $\mathcal{N} \cup \{K_2\}$,
      \item[\normalfont{\textbf{(T2)}}]\label{it:T2} there exists an algorithm $\mathrm{Alg}_{\mathcal{N}}$ that, for every $Q \in \mathcal{N} \cup \{K_2\}$ and $x : V(Q) \to \mathbb{Z}_{\geq 0}$, computes a minimum-size co-interval cover of $Q^x$, in time polynomial in $|E(Q^x)|$.
  \end{enumerate}
\end{definition}
Property~(T1) ensures that the block family is closed under vertex deletion.
Property (T2) guarantees an efficient algorithm on the pendant augmentations of graphs in $\mathcal{N} \cup \{K_2\}$.
As a special case where no pendant vertices are attached, this includes computing the co-boxicity of $Q$ itself.
Surprisingly, this local information about each block is enough to determine the co-boxicity of the entire graph.
\begin{theorem}
\label{thm:nice-block-algorithm}
Let $\mathcal{N}$ be a tractable family of $2$-connected graphs, and let $q$ be a polynomial bounding the running time of $\mathrm{Alg}_{\mathcal{N}}$.
If each block of a graph $G$ belongs to $\mathcal{N} \cup \{K_2\}$, then a minimum-size co-interval cover of $G$, and hence $\coboxicity(G)$, can be computed in time $O( |E(G)|^2 \cdot q(|E(G)|) )$.
\end{theorem}
A key idea underlying our proof of Theorem~\ref{thm:nice-block-algorithm} and establishing (T2) for certain graph classes, is that every co-interval subgraph of a graph $G$ is contained in an important type of subgraph, which we call a block ant.
A \emph{block ant} is a pendant augmentation of a block $Q$ in $G$ containing all edges incident to at most two vertices of $Q$.
After proving Theorem~\ref{thm:nice-block-algorithm}, we show that certain families of $2$-connected graphs are tractable, and then apply Theorem~\ref{thm:nice-block-algorithm} to each to expand the list of efficiently computable classes.
We start by showing that the empty family $\mathcal{N} = \emptyset$ is tractable.
Then, as blocks with more than two vertices contain a cycle, any tree with at least two vertices has blocks only in $\{K_2\}$, and so Theorem~\ref{thm:nice-block-algorithm} yields the following.
\begin{corollary}\label{cor:tree-coboxicity}
    Computing the co-boxicity of trees is in $\mathbf{P}$.
\end{corollary}
We extend Corollary~\ref{cor:tree-coboxicity} in two directions, by relaxing the order or the structure of the blocks.
First, we allow blocks of arbitrary structure but bounded order.
For every $k \geq 3$, the $2$-connected graphs with at least three and at most $k$ vertices form a tractable family.
The resulting running time is of the form $f(k) \cdot |E(G)|^2 + |E(G)|^3$, so the algorithm is fixed-parameter tractable when parameterized by $k$.
\begin{corollary}\label{cor:bounded-block-coboxicity}
    Let $k \geq 3$ and let $G$ be a graph whose blocks each have at most $k$ vertices.
    Then $\coboxicity(G)$ can be computed in time $O\big(k^{k^2 + 6k + 7} \cdot |E(G)|^2 + |E(G)|^3\big)$.
\end{corollary}
Second, we allow blocks of arbitrary order but impose block structure.
We show that $\mathcal{N} = \{K_k: k \geq 3\}$ is tractable.
This gives the following result for \emph{block graphs}, the class of graphs in which all blocks are from the set $\{K_k: k \geq 2\}$.
\begin{corollary}\label{cor:block-graph-coboxicity}
    Computing the co-boxicity of block graphs is in $\mathbf{P}$.
\end{corollary}
Block graphs are the graphs whose blocks are as dense as possible.
At the opposite extreme are cactus graphs, whose blocks are as sparse as possible.
A \emph{cactus graph} is a graph in which every block is a cycle or an edge.
We show that $\mathcal{N} = \{C_k: k \geq 3\}$ is tractable, and so we obtain the following.
\begin{corollary}\label{cor:cactus-coboxicity}
    Computing the co-boxicity of cactus graphs is in $\mathbf{P}$.
\end{corollary}
Tractability is preserved under unions, provided that one can efficiently recognize which family a block belongs to.
Indeed, if $\mathcal{N}_1$ and $\mathcal{N}_2$ are tractable and membership can be decided in polynomial time, then $\mathcal{N}_1 \cup \mathcal{N}_2$ is tractable: (T1) is inherited from the two families, and for (T2) it suffices to test whether $Q \in \mathcal{N}_1$ or $Q \in \mathcal{N}_2$  and run $\mathrm{Alg}_{\mathcal{N}_1}$ or $\mathrm{Alg}_{\mathcal{N}_2}$ accordingly.
As complete graphs, cycles, and graphs of bounded order are easy to recognize, combining the families above yields the following.
\begin{corollary}\label{cor:combined-coboxicity}
    Let $k \geq 3$.
    Computing the co-boxicity of graphs in which every block is a complete graph, a cycle, or has at most $k$ vertices is in $\mathbf{P}$.
\end{corollary}

The paper is organized as follows.
Section~\ref{section-prelims} reviews the background on co-interval graphs and co-boxicity.
Section~\ref{section-ECSs} establishes relevant block structure and shows key properties of maximal co-interval subgraphs related to blocks.
Section~\ref{sec:proof-of-main-theorem} proves Theorem~\ref{thm:nice-block-algorithm}, via a recursive algorithm. 
Section~\ref{sec:applications} proves that the four families above are tractable and derives Corollaries~\ref{cor:tree-coboxicity} to~\ref{cor:cactus-coboxicity}.
Section~\ref{section-discussion} concludes with open problems suggested by our method.
%

\section{Preliminaries}
\label{section-prelims}
In this section, we define basic graph terminology and preliminary theory on co-interval graph coverings used in this paper.

\subsection{Basic graph notation and terminology}
We use standard graph notation throughout (see~\cite{2003_Schrijver} for background and terminology).
All graphs in this paper are simple.
We briefly recall notation that will be used frequently.
Let $G = (V, E)$ be a graph.
For a vertex $v \in V$, the sets $\delta_G(v) \subseteq E$ and $N_G(v) \subseteq V$ are defined to be the edges incident to $v$ and vertices adjacent to $v$, respectively.
The \emph{closed neighbourhood} of $v$ is defined as $N_G[v] = N_G(v) \cup \{v\}$.
If $U \subseteq V$  is a set of vertices to be removed from $G$, then the induced subgraph $G[V \setminus U]$ is denoted by $G \setminus U$.
If $U = \{v\}$, then we write $G-v$ and $G\setminus U$ interchangeably.
Given two graphs $G_1 = (V_1, E_1)$ and $G_2 = (V_2, E_2)$, the graph $(V_1 \cup V_2, E_1 \cup E_2)$ is the \emph{union} of $G_1$ and $G_2$.
%

\subsection{Co-interval graphs and co-boxicity}\label{sec:subsec:co-intervalAndCoboxicity}

A graph $H = (U,F)$ is a \emph{co-interval graph} if there exists a function $I$ mapping $U$ to closed intervals in the real line such that $uv \in F$ if and only if $I(u) \cap I(v) = \emptyset$.
The function $I$ is called a \emph{co-interval representation} of $H$.
A class of graphs is \emph{hereditary} if all induced subgraphs remain in the class.
Given a co-interval graph $H$ with co-interval representation $I$, the subgraph induced by $S \subseteq V(H)$ is shown to be co-interval by restricting $I$ to the intervals corresponding to $S$.
So, the class of co-interval graphs is hereditary.
%

Co-interval graphs admit a forbidden induced subgraph characterization due to Lekkerkerker and Boland~\cite{1962_Lekkerkerker}.
Among the several forbidden induced subgraphs of co-interval graphs, only two play a central role in our analysis of block-structured graphs: $2K_2$, and $K_3$ with a single pendant attached to each vertex, denoted $\Gnet$. 
Moreover, in Section~\ref{sec:applications}, we use the fact that co-interval graphs induce no $C_5$ to prove Corollary~\ref{cor:cactus-coboxicity}.
\begin{lemma}[\cite{1962_Lekkerkerker}]\label{lem:co-interval-forbidden}
A co-interval graph contains no induced $2K_2$, $\Gnet$, or $C_5$.
\end{lemma}
\begin{definition}
A family $\cC = \{G_1, G_2, \ldots, G_t\}$ is a \textbf{co-interval cover} of a graph $G$ if each $G_i$ is a co-interval subgraph of $G$ and $\bigcup_{i = 1}^t E(G_i) = E(G)$.
\end{definition}
A co-interval cover reformulates an intersecting box representation of $G$ in the following sense.
Suppose $\coboxicity(G) = t$ so that there is a set of axis-parallel boxes $\{B(u) \subset \R^t: u \in V\}$ such that $uv \in E$ if and only if $B(u) \cap B(v) = \emptyset$.
Then, for each $1 \leq i \leq t$, the projection of these boxes onto the $i$-th axis forms a co-interval representation of a graph $G_i$ on $V$.
Since two axis-parallel boxes are disjoint if and only if their projections onto some axis are disjoint, each $G_i$ is a subgraph of $G$ and $\bigcup_{i = 1}^t E(G_i) = E$, so $\{G_1, G_2, \ldots, G_t\}$ is a co-interval cover of $G$.
Conversely, suppose $\{G_1, G_2, \ldots, G_t\}$ is a co-interval cover of $G$. 
As adding an interval to a co-interval representation that intersects all other intervals corresponds to adding an isolated vertex, we may assume that each $G_i$ spans $G$.
Then for each $u \in V$, the Cartesian product of the intervals representing $u$ in each element of the cover produces an axis-parallel box $B(u)$ in $\R^t$ such that for all $u,v \in V$, $uv \in E$ if and only if $B(u) \cap B(v) = \emptyset$.
So, when $t$ is smallest, $\coboxicity(G) = t$.
This reformulation proves the following equation.

\begin{lemma}[Cozzens and Roberts \cite{1983_Cozzens}]\label{lemma-boxicityIntervalOrderCover}
    Let $G$ be a graph. Then $\coboxicity(G)$ equals the minimum size of a co-interval cover of $G$.
\end{lemma}

Since co-interval graphs are hereditary, a co-interval cover of $G$ can be restricted to a co-interval cover of any induced subgraph $H$.
\begin{observation}
\label{obs:coboxicity-monotone}
    An induced subgraph $H$ of $G$ satisfies $\coboxicity(H) \leq \coboxicity(G)$.
\end{observation}
Given Lemma~\ref{lem:co-interval-forbidden}, as edges in distinct connected components of $G$ induce $2K_2$, each graph in a co-interval cover of $G$ has edges in at most one component.
By Lemma~\ref{lemma-boxicityIntervalOrderCover}, $\coboxicity(G)$ sums the co-boxicities of $G$'s connected components.
\begin{lemma}[Trotter~\cite{1979_Trotter}]\label{lem:box_and_disjoint_union}
    Let $H_1, H_2, \ldots, H_t$ be graphs with pairwise disjoint vertex sets and let $G = \bigcup_{i = 1}^t H_i$.
    Then $\coboxicity(G) = \sum_{i=1}^t \coboxicity(H_i)$.
\end{lemma}

\section{Block structure and maximal co-interval subgraphs}\label{section-ECSs}
By Lemma~\ref{lemma-boxicityIntervalOrderCover}, the co-boxicity of a graph is the size of a minimum-size co-interval cover.
Moreover, we may assume each co-interval graph is edge-wise maximal.
The first step in our covering approach is to establish that every maximal co-interval subgraph is contained in what we call a \emph{block ant}: a pendant augmentation of a block by $x$, where $x$ maps at most two vertices to positive integers.
Given a tractable class satisfying Definition~\ref{def:tractable}, optimal coverings within block ants can be obtained in polynomial time.
Hence, our focus is on describing how these coverings can be combined to efficiently cover the whole graph.
We begin by setting up important terminology regarding the blocks in a graph.

Any two blocks of a graph $G$ share at most one vertex, and the shared vertex must be a cut-vertex.
Two distinct blocks are \emph{neighbours} if they share a common cut-vertex.
We classify a block by the number of cut-vertices it contains: 
A block is \emph{isolated} if it contains no cut-vertex, a \emph{leaf} if it contains exactly one cut-vertex, and \emph{internal} otherwise.
An \emph{edge block} is a block containing exactly $2$ vertices.
An \emph{edge-leaf} block is a block that is both an edge block and a leaf block.
A block is \emph{substantial} if it is not an edge-leaf block.

\subsection{Block structure in graphs}\label{sec:BlockGraphs}
The blocks and cut-vertices of a connected graph exhibit a tree incidence structure, called the \emph{block-cut tree}~\cite{1966_Harary}.
The following observation about the blocks of a connected graph mirrors the elementary fact that every tree with more than one vertex has a leaf vertex.
\begin{observation}\label{obs:clique_or_2_leaves}
    If $G$ is a connected graph that is not an isolated block, then $G$ contains a leaf block.
\end{observation}
\begin{definition}\label{def:anchor}
  Given a substantial block $B$ of $G$, a cut-vertex $v \in V(B)$ is an \textbf{anchor} of $B$ if $v$ is shared with another substantial block.
\end{definition}

When $G$ contains an anchor, our approach identifies a substantial block in $G$ with exactly one anchor.
\begin{definition}\label{def:peripheral-block}
  A substantial block $Q$ of $G$ is a \textbf{peripheral block} if it has exactly one anchor.
\end{definition}
\begin{example}\label{example:block-graph-with-near-leaf-blocks}
    Figure~\ref{fig:BC_tree} shows a graph $G$ containing a block $Q$ as a peripheral block with anchor $u$.
    After deleting $V(Q) \setminus \{u\}$, the vertex $u$ is not an anchor anymore (nor a cut-vertex) and the substantial block $Q'$ becomes peripheral with anchor $u'$.
    \begin{figure}
        \centering
        \includegraphics[scale=0.85]{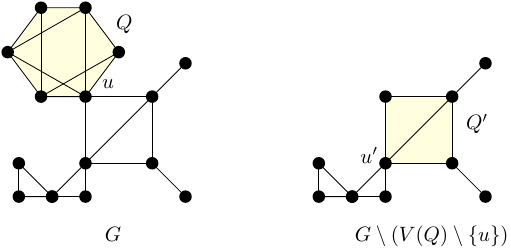}
        \caption{Peripheral blocks.}
        \label{fig:BC_tree}
    \end{figure}
\end{example}
\begin{lemma}\label{lem:peripheral-block-exists}
  Let $G$ be a connected graph.
  If $G$ has an anchor, then $G$ contains a peripheral block.
\end{lemma}
\begin{proof}
    Let $v$ be an anchor in $G$.
    Let $G'$ be obtained from $G$ by deleting all pendants of $G$.
    This removes all edge-leaf blocks from $G$ so that the blocks of $G'$ are the substantial blocks in $G$.
    Pendants are not cut-vertices, so $G'$ is connected.
    Since two substantial blocks share $v$ in $G$, $v$ is a cut-vertex in $G'$, and so $G'$ is not an isolated block. 
    By Observation~\ref{obs:clique_or_2_leaves}, in $G'$, there is a leaf block $B$ with cut-vertex $u$. 
    The vertex $u$ is shared by $B$ and some other block in $G'$, both of which are substantial in $G$, and so $u$ is an anchor in $G$.
    Since $B$ is a leaf block in $G'$, $u$ is the unique anchor of $B$ in $G$. 
    Thus, $B$ is a peripheral block of $G$.
\end{proof}
%

\subsection{Block ants and maximal co-interval subgraphs}\label{sec:subsec:ants}
\begin{figure}[htbp]
    \centering
    \includegraphics[scale=0.85]{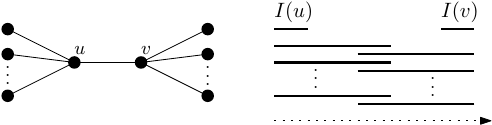}
    \caption{Co-interval representation $I$ of a $(u,v)$-ant.}
    \label{fig:ant-co-interval-representation}
\end{figure}

Ants are a class of co-interval graphs originally used in the context of co-boxicity in~\cite{1983_Cozzens} to construct co-interval covers.
Given a specified edge $uv$ of a graph $G$, the \emph{$(u,v)$-ant} subgraph of $G$ is $(N_G[u] \cup N_G[v], \delta_G(u) \cup \delta_G(v))$.
Figure~\ref{fig:ant-co-interval-representation} shows a co-interval representation of an arbitrary ant.
We generalize ants by replacing the specified edge with a block $Q$ of $G$.
\begin{definition}\label{definition-crab}
    Let $G$ be a graph, $Q$ a block of $G$, and $u,v \in V(Q)$.
    Then the subgraph
    \[Q_{u,v}(G) = \left( V(Q) \cup N_G(u) \cup N_G(v),\, E(Q) \cup \delta_G(u) \cup \delta_G(v) \right)\]
    is a \textbf{block ant} of $G$.
    We say that $Q_{u,v}(G)$ {\bf extends} $Q$ and call $u,v$ its {\bf hubs}.
    If $G$ is clear from the context, we write $Q_{u,v}$ instead of $Q_{u,v}(G)$.
\end{definition}
Note that if $u \neq v$, then, since $u, v \in V(Q)$, any common neighbour of $u$ and $v$ lies in $V(Q)$.
By Lemma~\ref{lem:co-interval-forbidden}, a co-interval subgraph of $G$ has no induced $2K_2$, and so all of its edges are contained within one connected component of $G$.
So, every co-interval subgraph of $G$ consists of a single connected co-interval subgraph with possibly some isolated vertices.
Since we are covering edges, it thus suffices to locate connected co-interval subgraphs of $G$.
A co-interval subgraph $H$ of $G$ is \emph{maximal} if there does not exist a co-interval subgraph $K$ of $G$ satisfying $E(H) \subsetneq E(K)$.
We show three facts about maximal co-interval subgraphs necessary to prove Theorem~\ref{thm:nice-block-algorithm}.

\begin{figure}[htbp]
    \centering
    \includegraphics[scale=0.85]{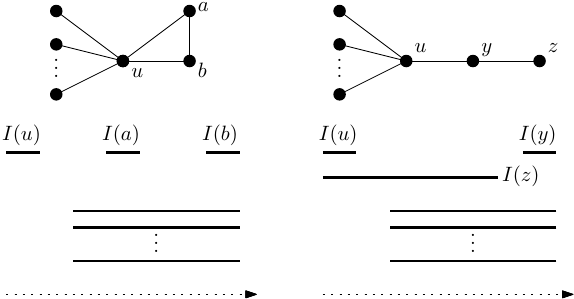}
    \caption{Cases~\textit{(i)} and~\textit{(ii)} of $H'$ in the proof of Proposition~\ref{prop:maximal-star-or-more} with co-interval representation $I$.}
    \label{fig:Aug_Stars}
\end{figure}

\begin{proposition}\label{prop:maximal-star-or-more}
Let $G$ be a connected graph that is not a star.
Then no maximal co-interval subgraph of $G$ is a star.
\end{proposition}
\begin{proof}
Suppose for a contradiction that $G$ has a maximal co-interval subgraph $H$ that is a star with centre $u$.
As any subgraph of a star is co-interval, if some edge $e \in \delta_G(u)$ satisfies $e \notin E(H)$, then $H+e$ contradicts the maximality of $H$.
Thus, $E(H) = \delta_G(u)$.
Since $G$ is connected and not a star, either \textit{(i)} two vertices $a,b \in N_G(u)$ are adjacent in $G$, or \textit{(ii)} for some two vertices $y$ and $z$, $G[\{u,y,z\}]$ is a path.
If~\textit{(i)}, set $H' := H \cup (\{a,b\},\{ab\})$.
If~\textit{(ii)}, set $H' := H \cup (\{y,z\},\{yz\})$.
Then Figure~\ref{fig:Aug_Stars} shows a co-interval representation $I$ of $H'$, implying that $H'$ is a co-interval graph.
Since $ E(H) \subsetneq E(H') \subseteq E(G)$, $H'$ contradicts the maximality of $H$.
\end{proof}

We show in Lemma~\ref{lem:ECS_are_subgraph_of_big_ant} that every connected co-interval subgraph $H$ of $G$ is contained within some block ant of $G$.
First, we prove the following supporting lemma.
\begin{lemma}\label{lem:degree_two_in_Q}
    Let $H$ be a connected co-interval subgraph of a graph $G$.
    Then all vertices of degree at least two in $H$ belong to a common block of $G$.
\end{lemma}
\begin{proof}
    Suppose for a contradiction that $u$ and $v$ are vertices not in the same block of $G$ with degree at least $2$ in $H$.
    As no block contains both $u$ and $v$, $uv \notin E(G)$.
    Moreover, since $H$ is connected, there is a cut vertex $c$ separating them.
    Since $u$ and $v$ have degree at least $2$ in $H$, there exist $n_u \in N_H(u) \setminus \{c\}$ and $n_v \in N_H(v) \setminus \{c\}$.
    As neighbours of $u$ and $v$ other than $c$, $n_u$ and $n_v$ lie in different blocks, so $n_u \neq n_v$ and $n_u n_v, u n_v, v n_u \notin E(G)$.
    Since $uv \notin E(G)$, the set $\{u,n_u,v,n_v\}$ induces $2K_2$ in $H$, contradicting Lemma~\ref{lem:co-interval-forbidden}.
\end{proof}
Lemma~\ref{lem:ECS_are_subgraph_of_big_ant} shows that every connected co-interval subgraph lies inside some block ant.
\begin{lemma}\label{lem:ECS_are_subgraph_of_big_ant}
    Let $H$ be a connected co-interval subgraph of a graph $G$ such that $H$ has at least one edge.
    Then $H$ is a subgraph of a block ant of $G$.
\end{lemma}
\begin{proof}
Suppose $H$ has no vertex of degree at least two.
Since $H$ is connected, it is $K_2$, and so it is in exactly one block.
Otherwise, let $Q$ be the block containing all vertices of degree at least two in $H$, given by~Lemma~\ref{lem:degree_two_in_Q}.
We show the following:
\begin{equation}\label{eq:two_active_vertices}
    \text{At most two vertices of $Q$ are incident to edges in $E(H) \setminus E(Q)$.}
\end{equation}
Suppose for a contradiction that $Q$ contains three distinct vertices $u_1, u_2, u_3$, each incident to an edge in $E(H) \setminus E(Q)$, with respective neighbours $v_1, v_2, v_3 \notin V(Q)$.
By~Lemma~\ref{lem:degree_two_in_Q}, $v_i$ has degree $1$ in $H$, with unique edge $v_iu_i$ in $H$.
If two distinct $i,j \in \{1,2,3\}$ satisfy $u_iu_j \notin E(H)$, then $\{u_i,v_i,u_j,v_j\}$ induces $2K_2$ in $H$.
Otherwise, $\{u_1,u_2,u_3,v_1,v_2,v_3\}$ induces $\Gnet$ in $H$.
Both cases contradict Lemma~\ref{lem:co-interval-forbidden}.

\smallskip
Using~\eqref{eq:two_active_vertices}, fix at most two vertices $u, v \in V(Q)$ such that every vertex of $Q$ incident to an edge of $E(H) \setminus E(Q)$ is $u$ or $v$.
By~Lemma~\ref{lem:degree_two_in_Q} and the connectivity of $H$, each vertex in $V(H) \setminus V(Q)$ is
a leaf of $H$ with unique neighbour in $V(Q)$.
Thus, $E(H) \subseteq E(Q) \cup \delta_G(u) \cup \delta_G(v)$ and $H$ is a
subgraph of $Q_{u,v}$.
\end{proof}

Second, we show that every maximal connected co-interval subgraph lies inside a block ant extending a block $Q$, where $Q$ is unique and substantial.
\begin{lemma}\label{lem:unique-ant}
Let $G$ be a connected graph that is not a star and $H$ be a maximal connected co-interval subgraph of $G$.
Then there is a unique block $Q$ such that $H$ is a subgraph of a block ant that extends $Q$.
Moreover, $Q$ is a substantial block.
\end{lemma}
\begin{proof}
By Lemma~\ref{lem:ECS_are_subgraph_of_big_ant}, there is a block $Q$ in $G$ such that $H$ is a subgraph of some block ant that extends $Q$.
For uniqueness, suppose for a contradiction that $H$ is contained in $Q_{u,w}$ and $Q'_{u',w'}$ for distinct blocks $Q \neq Q'$.
If $H$ has an edge $e \in E(Q)$, then since blocks partition the edges of $G$, $e \notin E(Q')$, so $e$ must lie in $\delta_G(u') \cup \delta_G(w')$.
Thus, one endpoint of $e$ is $u'$ or $w'$.
That endpoint also lies in $V(Q)$, and since $Q, Q'$ share at most one vertex, this forces $u'$ (or $w'$) to equal the unique shared vertex $c^*$ in $V(Q) \cap V(Q')$.
Symmetrically, if $H$ has an edge in $E(Q')$, the same argument forces $u$ (or $w$) to equal $c^{*}$.
Thus, every edge of $H$ lying in $E(Q)$ or $E(Q')$ is then incident to $c^{*}$.
Moreover, every remaining edge of $H$ must be incident to a vertex of both ants.
If such an edge avoided $c^{*}$, it would have to join a vertex of $V(Q)\setminus\{c^{*}\}$ to a vertex of $V(Q')\setminus\{c^{*}\}$, which is impossible since $Q$ and $Q'$ are distinct blocks.
So $H$ is a subgraph of the star centred at $c^{*}$.
By maximality of $H$ and since stars are co-interval, $H$ is itself the star centred at $c^{*}$, contradicting Proposition~\ref{prop:maximal-star-or-more}.

Finally, suppose $Q$ is an edge-leaf block, i.e. $Q = \{c, \ell\}$ where $c$ is a cut-vertex and $\ell$ is a leaf.
Thus, for any $u, w \in V(Q)$, the edge set of $Q_{u,w}$ is contained in $\delta_G(c)$.
Once more, by maximality, $H$ is itself the star centred at $c$, contradicting Proposition~\ref{prop:maximal-star-or-more}.
\end{proof}

Third, having pinned down the block $Q$, we show that maximality forces $H$ to contain all edges incident to any hub it reaches outside of $Q$.
The proof below uses the \emph{Helly property} for intervals: any family of pairwise intersecting intervals on $\mathbb{R}$ has a common point~\cite{2004_Golumbic}.
\begin{lemma}
\label{lem:hub-saturation}
Let $H$ be a maximal connected co-interval subgraph of a graph $G$. 
Suppose $H$ is contained in a block ant with hub $u$ that extends a block $Q$ of $G$.
If there exists $v \in V(G) \setminus V(Q)$ such that $uv \in E(H)$, then $\delta_G(u) \subseteq E(H)$.
\end{lemma}
\begin{proof}
    Suppose for a contradiction that there exists $y \in N_G(u)$ with $uy \in E(G) \setminus E(H)$.
    We show that $H' := (V(H) \cup \{y\}, E(H) \cup \{uy\})$ is a co-interval graph, contradicting the maximality of $H$.
    Let $\{I(z): z \in V(H)\}$ be a co-interval representation of $H$.
    Since $Q$ is a block, every vertex in $V(G) \setminus V(Q)$ has at most one neighbour in $V(Q)$.
    Thus, $u$ is the only neighbour of $v$ in $H$.
    We have
    \begin{equation}\label{eq:x-saturated}
        I(u) \cap I(v) = \emptyset \quad\text{and}\quad I(v) \cap I(z) \neq \emptyset \text{ for every } z \in V(H) \setminus \{u\}.
    \end{equation}
    If $y \notin V(H)$, then, by \eqref{eq:x-saturated}, setting $I(y) := I(v)$ extends $I$ to a co-interval representation of $H'$, since $u$ is then the only neighbour of $y$ in $H'$.
    Otherwise, assume $y \in V(H)$.
    Replace $I(y)$ with $I'(y) := I(y) \cap I(v)$, leaving $I(z)$ unchanged for $z \neq y$.
    We verify that $I$ is a co-interval representation of $H'$.
    Since we only change the interval corresponding to $y$, it is sufficient to check the following:
    \begin{equation*}
        \text{for every $z \in V(H)$, $zy \in E(H')$ if and only if $I(z) \cap I'(y) = \emptyset$.}
    \end{equation*}
    Recall that $E(H') = E(H) \cup \{uy\}$.
    For $uy \in E(H')$, we have $I'(y) \cap I(u) \subseteq I(x) \cap I(u) = \emptyset$. 
    For $zy \in E(H)$, we have $I'(y) \cap I(z) \subseteq I(y) \cap I(z) = \emptyset$. 
    Now let $zy \notin E(H')$.
    Since $uy \in E(H')$, $z\neq u$.
    Hence, $zy \notin E(H)$ and so $I(z) \cap I(y) \neq \emptyset$.
    Moreover, \eqref{eq:x-saturated} yields $I(z) \cap I(x) \neq \emptyset$ and $I(y) \cap I(v) \neq \emptyset$.
    As closed intervals satisfy the Helly property and $I(z), I(v), I(y)$ are pairwise intersecting, $I(z) \cap I(y) \cap I(v) \neq \emptyset$, and in particular $I(z) \cap I'(y) \neq \emptyset$.
\end{proof}

\section{Proof of Theorem~\ref{thm:nice-block-algorithm}}
\label{sec:proof-of-main-theorem}
We prove Theorem~\ref{thm:nice-block-algorithm} via a recursive procedure, which we call $\mathrm{Cov}$ and describe in Algorithm~\ref{alg:coboxicity}.
By Definition~\ref{def:tractable}, $\mathcal{N}$ comes equipped with an algorithm $\mathrm{Alg}_{\mathcal{N}}$.
At a high level, $\mathrm{Cov}$ considers a connected component $K$ of its input graph.
If $K$ is itself a pendant augmentation, $\mathrm{Alg}_{\mathcal{N}}$ applies to it directly.
Otherwise, $K$ contains a peripheral block $B$ intersecting the rest of the graph at its anchor $v$.
Probing the neighbourhood of $v$, the covering problem is split into two nearly-disjoint pieces, one involving $B$ and the other involving the remainder of the graph.
The key subtlety is deciding exactly how much of $v$'s neighbourhood to fold into which of the two pieces.
The source of efficiency of the algorithm arises from the fact that it is optimal to cover $K$ by combining a cover of one of just three local candidates built around $v$ and $B$ with an optimal cover of the remainder of the graph.
Moreover, which of these three candidates is the right one is decided efficiently by comparing the sizes of the three local covers, rather than by exploring the potentially much larger remainder of the graph.
%
\renewcommand{\algorithmicrequire}{\textbf{Input:}}     
\renewcommand{\algorithmicensure}{\textbf{Output:}}     
\begin{algorithm}[ht]
\caption{$\mathrm{Cov}(\Gamma)$ -- a minimum-size co-interval cover of $\Gamma$}
\label{alg:coboxicity}
\begin{algorithmic}
\setlength{\itemsep}{2pt}
\Require A graph $\Gamma$ such that every block of $\Gamma$ belongs to $\mathcal{N} \cup \{K_2\}$.
\Ensure A minimum-size co-interval cover of $\Gamma$.
\If{$E(\Gamma) = \emptyset$}
    \State \Return $\emptyset$
\EndIf
\State Let $K$ be a connected component of $\Gamma$ containing an edge.
\If{$K$ has no anchor} \Comment{Case 1}
    \State Write $K = B^{x}$ for a block $B$ of $K$, where $x_v = |N_K(v)\setminus V(B)|$ for each $v \in V(B)$.
    \State \Return $\mathrm{Alg}_{\mathcal{N}}(B,x) \cup \mathrm{Cov}(\Gamma \setminus V(K))$
\EndIf
\State Let $B$ be a peripheral block of $K$ with anchor $v$. \Comment{Case 2}
\State Let $H^-$ be the component of $K-v$ containing $V(B)\setminus\{v\}$.
\State Set  $H := B^{x}$, where $x_v = 0$ and $x_u = |N_{H^-}(u)\setminus V(B)|$ for $u \in V(B) \setminus \{v\}$.
\State Set $H^+ := B^{y}$, where $y_v = |N_K(v)\setminus V(H)|$ and $y_u = x_u$ for $u \in V(B) \setminus \{v\}$.
\State Compute $\cC_{H^-} := \mathrm{Cov}(H^-)$, $\cC_H := \mathrm{Alg}_{\mathcal{N}}(B,x)$ and $\cC_{H^+} := \mathrm{Alg}_{\mathcal{N}}(B,y)$.
\If{$|\cC_{H^-}| = |\cC_{H^+}|$} \Comment{Case 2(a)}
    \State \Return $\cC_{H^+} \cup \mathrm{Cov}(\Gamma \setminus V(H))$
\EndIf
\If{$|\cC_{H^-}| = |\cC_H|$} \Comment{Case 2(b)}
    \State \Return $\cC_H \cup \mathrm{Cov}(\Gamma \setminus V(H^-))$
\EndIf
\State Let $\Gamma^*$ be $\Gamma$ with $E(H^-)$ and $V(H^-)\setminus N_K[v]$ removed. \Comment{Case 2(c)}
\State \Return $\cC_{H^-} \cup \mathrm{Cov}(\Gamma^*)$
\end{algorithmic}
\end{algorithm}
%
\subsection{$\mathrm{Cov}$ is well-defined}
This subsection shows that every step of $\mathrm{Cov}$ 
is legitimate: whichever case applies, the graphs it passes to $\mathrm{Alg}_{\mathcal{N}}$ are pendant augmentations of a block, as required by (T2), and the graphs it recurses on satisfy the hypothesis of
$\mathrm{Cov}$.

%
In Case~1, $K$ has no anchor.
The next lemma shows that this forces $K$ to be a pendant augmentation of a block.
If there is an isomorphism from $Q^{x}$ to $G$ that is the identity on $V(Q)$, we write, with slight abuse of notation, $G = Q^{x}$.
\begin{lemma}\label{lem:no-anchor-block-plus-pendants}
    Let $G$ be a connected graph. 
    If $G$ has no anchor, then $G = Q^{x}$ for some block $Q$ of $G$, where $x$ maps each $v \in V(Q)$ to $|N_G(v) \setminus V(Q)|$, the number of block neighbours of $Q$ containing $v$.
\end{lemma}
\begin{proof}
    We identify a \emph{suitable} block $Q$ of $G$ for the claimed $x$.
    If $G$ is a block, then $Q := G$ is suitable.
    Suppose $G$ has multiple blocks. 
    Assume each block of $G$ is an edge-leaf. 
    As $G$ is connected, $G$ is a star centred at $u$. 
    Then for any $w \in V(G) \setminus \{u\}$, the block $Q := (\{u,w\}, \{uw\})$ is suitable.
    Suppose $G$ has a substantial block $B$. 
    As $G$ has no anchor, each block neighbour of $B$ is an edge-leaf. 
    So, for any $v \in V(B)$, the number of block neighbours of $B$ containing $v$ is $|N_G(v) \setminus V(B)|$.
    As $G$ is connected, $G = B^x$ and hence $Q := B$ is suitable.
\end{proof}

In Case~2, $K$ has an anchor, so by Lemma~\ref{lem:peripheral-block-exists}, $K$ contains a peripheral block $B$ with anchor $v$.
The next lemma shows that $H$, as constructed from $B$ and $v$, is itself a pendant augmentation of $B$.
Note that $H^{+}$ is, by construction, a pendant augmentation of $B$ as well.
\begin{lemma}\label{lem:peripheral-block-structure}
	Let $B$ be a peripheral block of $G$ with anchor $v$.
	Then $G - v$ has a unique connected component $H^-$ containing $V(B) \setminus \{v\}$ such that
	$
	G[V(H^-) \cup \{v\}] = B^{x},
	$
	where $x_v := 0$ and $x_u := |N_G(u) \setminus V(B)|$ for all $u \in V(B) \setminus \{v\}$.
\end{lemma}
\begin{proof}
    As $B$ is $2$-connected, $B-v$ is connected, so $V(B) \setminus \{v\}$ is in a connected component $H^-$ of $G - v$.
    As $B$ is peripheral with anchor $v$, each block neighbour of $B$ not sharing $v$ is an edge-leaf in $G$. 
    It follows that $V(H^-)$ is the union of the vertices in $V(B) \setminus \{v\}$ with their pendants.
    So, $G[V(H^-) \cup \{v\}] = B^x$.
\end{proof}

We show that every graph on which $\mathrm{Cov}$ recurses also satisfies the hypothesis of $\mathrm{Cov}$.
This is exactly where Property~(T1) of Definition~\ref{def:tractable} is needed.
Most of the recursive calls simply restrict to an induced subgraph of $\Gamma$, so we first record that this operation preserves the algorithm's hypothesis.
\begin{lemma}\label{lem:hereditary}
If every block of a graph $\Gamma$ belongs to $\mathcal{N} \cup \{K_2\}$, then for every $S \subseteq V(\Gamma)$, every block of $\Gamma[S]$ also belongs to $\mathcal{N} \cup \{K_2\}$.
\end{lemma}
\begin{proof}
    Since $\Gamma[S]$ is obtained from $\Gamma$ by deleting the vertices of $V(\Gamma) \setminus S$ one at a time, it suffices to prove the claim for $S = V(\Gamma) \setminus \{u\}$, where $u \in V(\Gamma)$.
    Let $Q'$ be a block of $\Gamma - u$.
    As $Q'$ is a $2$-connected subgraph of $\Gamma$, it is contained in a block $Q$ of $\Gamma$.
    If $u \notin V(Q)$, then $Q$ is a $2$-connected subgraph of $\Gamma - u$ which, by the maximality of $Q'$, implies $Q' = Q \in \mathcal{N} \cup \{K_2\}$.
    Suppose $u \in V(Q)$.
    Since $Q'$ is a block of $Q - u$, $Q'$ has an edge, and so $Q \neq K_2$.
    By (T1), $Q' \in \mathcal{N} \cup \{K_2\}$.
\end{proof}

\begin{lemma}\label{lem:iteration-invariant}
    If every block of $\Gamma$ belongs to $\mathcal{N} \cup \{K_2\}$, then the same holds for every graph on which $\mathrm{Cov}$ recurses while processing $\Gamma$.
\end{lemma}
\begin{proof}
It suffices to check the graphs on which $\mathrm{Cov}$ makes a recursive call directly while processing $\Gamma$.
Let $K$ be the component of $\Gamma$ processed by $\mathrm{Cov}$; since the rest of $\Gamma$ is left untouched by this recursive call, it suffices to check the graphs derived from $K$.
In Case~1, Case~2(a), Case~2(b), and for the call on $H^{-}$, the resulting graph is an induced subgraph of $\Gamma$, so Lemma~\ref{lem:hereditary} applies directly.
In Case~2(c), the recursive call is on $\Gamma^*$, obtained by removing the edges of $E(H^{-})$ and the vertices of $V(H^{-}) \setminus N_K[v]$.
Equivalently, $\Gamma^*$ is obtained from the induced subgraph $\Gamma \setminus V(H^{-})$, for which every block belongs to $\mathcal{N} \cup \{K_2\}$ by Lemma~\ref{lem:hereditary}, by reattaching, at $v$, one pendant vertex for each retained vertex of $N_K(v) \cap V(H^{-})$, since each such vertex loses every edge except the one to $v$ once $E(H^{-})$ is removed.
Attaching a pendant vertex at $v$ only introduces a new edge-leaf block, in $\{K_2\}$, without altering any other block of the graph, so this too preserves the hypothesis.
\end{proof}

\subsection{$\mathrm{Cov}$ is polynomial}
\label{subsec:alg-coboxicity-is-polynomial}
To compute the time complexity of Algorithm~\ref{alg:coboxicity}, we focus on the edges of $K$ that the algorithm processes.
By (T2), Case~1 takes a time polynomial in $|E(K)|$.
In Case~2, removing the cut-vertex $v$ splits $K - v$ into two parts: $H^-$ and the rest of $K$. 
Write $n := |E(K)|$ and $m := |E(H^-)|$.
\begin{lemma}\label{lem:progress}
In Case~2, we have $1 \le m < n$.
Moreover,
\[
|E(K\setminus V(H))| \ \le\ |E(K\setminus V(H^-))| \ \le\ |E(K)\setminus E(H^-)| \ =\ n - m.
\]
\end{lemma}

\begin{proof}
Observe that the edges of $K$ incident to $v$ split into those going to $H^-$ and those going to the rest of $K$. 
Together with the edges lying entirely inside each of the two parts, this gives a partition of $E(K)$ into four sets:
\begin{equation}\label{eq:partition}
    E(K) = E(H^-) \ \cup \ \delta_{H}(v) \ \cup \ \delta_{K \setminus V(H^-)}(v) \ \cup \ E(K\setminus V(H)).
\end{equation}
Here, $\delta_H(v)$ is the set of edges joining $v$ and $H^-$, and $\delta_{K\setminus V(H^-)}(v)$ the remaining edges incident to $v$.
Since $B$ is substantial, it is not an edge-leaf block: either~\textit{(i)} $|V(B)| \ge 3$, or~\textit{(ii)} $B = (\{u,v\},\{uv\})$, where $u$ is a cut-vertex in $K$.
If~\textit{(i)}, as $B$ is $2$-connected, $B - v$ contains an edge, and this edge lies in $H^-$.
If~\textit{(ii)}, since $u$ is a cut-vertex of $K$ other than $v$, it has a neighbour outside $B$, and since $u \in V(H^-)$ this neighbour lies in $H^-$ as well, again giving an edge of $H^-$.
So, \textit{(i)}~and~\textit{(ii)} each imply $m \ge 1$.
Since $v$ is a cut-vertex of $K$ joining $H$ with $K \setminus V(H^-)$, $|\delta_H(v)| \geq 1$ and $|\delta_{K \setminus V(H^-)}(v)| \geq 1$.
So by \eqref{eq:partition},
$$n = m + |\delta_H(v)| + |\delta_{K \setminus V(H^-)}(v)| + |E(K\setminus V(H))| \geq m + 1 + 1 + 0 > m.$$
The inequality chain follows from the observation that $E(K \setminus V(H^-))$ is equal to both $E(K \setminus V(H)) \cup \delta_{K \setminus V(H^-)}(v)$ and $(E(K) \setminus E(H^-)) \setminus \delta_H(v)$.
\end{proof}

Recall that~(T2) yields a polynomial $q$ such that $\mathrm{Alg}_{\mathcal{N}}$ runs in time at most $q(|E(Q^x)|)$ for every $Q \in \mathcal{N} \cup \{K_2\}$ and every $x : V(Q) \to \mathbb{Z}_{\geq 0}$.
\begin{proposition}\label{prop:polynomial}
Let $G$ be a valid input for $\mathrm{Cov}$.
Then $\mathrm{Cov}(G)$ is computed in time $O(|E(G)|^2 \cdot q(|E(G)|))$.
\end{proposition}
\begin{proof}
Let $n_0 := |E(G)|$, and consider the recursion tree $\mathcal{T}$ for the computation of $\mathrm{Cov}(G)$, where the children of a parent $\mathrm{Cov}(\Gamma)$ are determined by the cases in Algorithm~\ref{alg:coboxicity} applied to $\Gamma$ as follows.
\textit{Case~1:} one child for $\mathrm{Cov}(\Gamma \setminus V(K))$; \textit{Case~2:} a child for $\mathrm{Cov}(H^-)$, and a second child for \textit{(a):} $\mathrm{Cov}(\Gamma \setminus V(H))$, \textit{(b):} $\mathrm{Cov}(\Gamma \setminus V(H^-))$, or \textit{(c):} $\mathrm{Cov}(\Gamma^*)$.
Label each node by $N := |E(\Gamma)|$ for the graph $\Gamma$ it processes.

In Case~1, since $K$ contains at least one edge, its unique child has size $N - n \le N-1$.
By Lemma~\ref{lem:progress}, in Case~2 the two children have sizes $m$ and at most $N-m$, since $1 \le m < n$ and the second child inherits at most $n-m$ edges of $K$ together with the $N-n$ edges of $\Gamma$ outside $K$.
In both cases, the children's sizes sum to at most $N$, and each is strictly smaller than $N$.
The strict decrease bounds the depth of $\mathcal{T}$ by $n_0$, since labels are non-negative integers strictly decreasing along any root-to-leaf path, starting at $n_0$.
Note that depth and branching factor alone only bound $\mathcal{T}$'s size by $2^{n_0}$; the polynomial bound instead comes from the children's sizes summing to at most the parent's.
Applying this level by level, the sum of the labels of the nodes at any fixed depth is at most $n_0$, the label of the root.
Since every non-trivial node ($N \ge 1$) contributes at least $1$ to its level's sum, each level contains at most $n_0$ non-trivial nodes with at most $n_0$ levels.
Thus, $\mathcal{T}$ has at most $n_0^2$ non-trivial nodes, plus at most $2n_0^2$ further nodes with $N=0$ (their children), each costing $O(1)$.

It remains to bound the cost of a single node, excluding its children.
A connected component $K$ of $\Gamma$ can be found by a single breadth-first search from a non-isolated vertex, in time $O(|E(K)|) = O(n)$.
Given $K$, a linear-time algorithm~\cite{1971_Paton} computes its block-cut tree, from which a single further linear scan determines whether $K$ has an anchor and, if so, identifies a peripheral block $B$ with its anchor $v$.
Constructing $H, H^-, H^+$ from $B$ and $v$ takes further linear time, and each of the at most two calls to $\mathrm{Alg}_{\mathcal{N}}$ costs time polynomial in $|E(B^x)|$ (or $|E(B^y)|$), by Property~(T2).
Since $E(B)$ and the pendant edges counted by $x$ (or $y$) are disjoint sets of edges of $K$, we have $|E(B^x)|, |E(B^y)| \le n \le n_0$, so there is a single polynomial $q$, independent of the node, bounding all of the above together with the $O(1)$ comparisons deciding the branch and the assembly of the returned cover.
The total running time is therefore $O(n_0^2 \cdot q(n_0))$.
\end{proof}

\subsection{$\mathrm{Cov}$ returns a minimum co-interval cover}
\label{subsec:alg-coboxicity-succeeds}
Let $\Gamma$ be a valid input for $\mathrm{Cov}$.
We prove, by strong induction on $N := |E(\Gamma)|$, that $\mathrm{Cov}(\Gamma)$ is a minimum-size co-interval cover of $\Gamma$.
By Lemma~\ref{lemma-boxicityIntervalOrderCover}, this shows $\coboxicity(\Gamma) = |\mathrm{Cov}(\Gamma)|$.
For the base case $N = 0$, $\mathrm{Cov}(\Gamma) = \emptyset$ is a minimum-size co-interval cover since $\Gamma$ has no edges. 
For the inductive step, fix $\Gamma$ with $N \ge 1$, and assume that for every valid input $\Gamma'$ with fewer than $N$ edges, $\mathrm{Cov}(\Gamma')$ returns a minimum-size co-interval cover of $\Gamma'$.
We proceed in reference to each case described in Algorithm~\ref{alg:coboxicity} to compute $\mathrm{Cov}(\Gamma)$.

\smallskip

\noindent
\textbf{Case 1.}
By Lemma~\ref{lem:no-anchor-block-plus-pendants}, $K = B^{x}$ for the unique block $B$ of $K$.
By Property~(T2), $\mathrm{Alg}_{\mathcal{N}}(B,x)$ is a minimum-size co-interval cover of $K$. 
Since $K$ contains an edge, $\Gamma\setminus V(K)$ has fewer than $N$ edges, so by the inductive hypothesis, $\mathrm{Cov}(\Gamma\setminus V(K))$ is a minimum-size co-interval cover of $\Gamma\setminus V(K)$.
As $\Gamma = K \cup (\Gamma\setminus V(K))$ and $E(K) \cap E(\Gamma\setminus V(K)) = \emptyset$, the returned family $\mathrm{Alg}_{\mathcal{N}}(B,x) \cup \mathrm{Cov}(\Gamma\setminus V(K))$ is a minimum-size co-interval cover of $\Gamma$.

\smallskip

\noindent
\textbf{Case 2.}
Let $t := |\cC_H|$ and $t^+ := |\cC_{H^+}|$.
By (T2), $t = \coboxicity(H)$ and $t^+ = \coboxicity(H^+)$.
As $v$ is a cut-vertex in $K$, $\delta_K(v) \neq \emptyset$, giving $|E(H^-)| < |E(K)| \leq N$. 
Thus, the inductive hypothesis gives $t^{-} := |\cC_{H^{-}}| = \coboxicity(H^{-})$.
\begin{lemma}
\label{lem:t-sandwich}
$t^{-} \le t \le t^{+} \le t^{-} + 1$.
\end{lemma}
\begin{proof}
The graphs $H^-$ and $H$ are induced subgraphs of $H$ and $H^+$, respectively.
Thus, by Observation~\ref{obs:coboxicity-monotone}, $t^- \leq t \leq t^+$.
As stars are co-interval graphs, and $E(H^+) \setminus E(H^-)$ can be covered by a star centred at $v$, adding such a star to $\cC_{H^-}$ produces a cover of $H^+$.
Thus, $t^+ \leq t^- + 1$.
\end{proof}
While processing $\Gamma$, let $\mathcal{P}$ denote the family of co-interval subgraphs Algorithm~\ref{alg:coboxicity} adds to the cover, and let $K'$ denote the subgraph that takes the place of $K$ in the graph passed to the recursive call.
In particular,
\[
(\mathcal{P},\, K') = \begin{cases}
    \bigl(\cC_{H^{+}},\ K \setminus V(H)\bigr), & \text{Case 2(a);} \\[2pt]
    \bigl(\cC_{H},\ K \setminus V(H^{-})\bigr), & \text{Case 2(b);} \\[2pt]
    \Bigl(\cC_{H^{-}},\ \bigl(V(K)\setminus(V(H^{-})\setminus N_K[v]),\ E(K)\setminus E(H^{-})\bigr)\Bigr), & \text{Case 2(c).}
\end{cases}
\]
Note that the recursive call is made on the graph $\Gamma' := K' \cup (\Gamma \setminus V(K))$.
Thus, Algorithm~\ref{alg:coboxicity} returns $\mathcal{P} \cup \mathrm{Cov}(\Gamma')$.
\begin{lemma}
\label{lem:iter-covers}
$\mathcal{P} \cup \mathrm{Cov}(\Gamma')$ is a co-interval cover of $\Gamma$.
\end{lemma}
\begin{proof}
By the inductive hypothesis, $\mathrm{Cov}(\Gamma')$ is a co-interval cover of $\Gamma' = K' \cup (\Gamma\setminus V(K))$. 
So, $\mathrm{Cov}(\Gamma')$ covers $E(K') \cup (E(\Gamma)\setminus E(K))$.
It is sufficient to show that $\mathcal{P}$ covers the remaining edges $E(K) \setminus E(K')$.
In Case~2(a), $E(K) \setminus E(K')$ is exactly the set of edges of $K$ incident to $V(H)$.
By Lemma~\ref{lem:peripheral-block-structure}, $H^{-}$ is the unique connected component of $K - v$ containing $V(B)\setminus \{v\}$.
So, every edge in $E(K) \setminus E(K')$ either lies in $E(H)$ or is incident to $v$ and a vertex outside $V(H)$.
That is, $E(K)\setminus E(K') = E(H) \cup \delta_{K\setminus V(H)}(v) = E(H^{+})$.
As $\mathcal{P} = \cC_{H^+}$ by definition, it covers $E(H^+)$, and hence $\mathcal{P}$ covers $E(K) \setminus E(K')$.
In Case~2(b), $E(K)\setminus E(K')$ consists of the edges of $K$ incident to $V(H^{-})$ which are covered by $\mathcal{P} = \cC_H$.
In Case~2(c), $E(K)\setminus E(K') = E(H^{-})$, covered by $\mathcal{P} = \cC_{H^-}$.
\end{proof}
We now prove optimality. 
As $\mathrm{Cov}(\Gamma')$ is optimal by the inductive hypothesis, we must show that $\mathcal{P}$ accounts for exactly the drop in co-boxicity from $K$ to $K'$.
\begin{lemma}
\label{lem:iter-decreases}
$\coboxicity(K') +  |\mathcal{P}| \le \coboxicity(K)$.
\end{lemma}
\begin{proof}
Let $\mathcal{C}$ be some minimum-size co-interval cover of $K$ consisting of maximal connected co-interval subgraphs.
Then $|\mathcal{C}| = \coboxicity(K)$.
Define
\[
\mathcal{C}(H^{-}) := \{\, C \in \mathcal{C} : E(C) \cap E(H^{-}) \neq \emptyset \,\}.
\]
We prove three claims about $\mathcal{C}$ and $\mathcal{C}(H^{-})$, and use them to produce, in each case, a co-interval cover $\mathcal{C}'$ of $K$ with $|\mathcal{C}'| = |\mathcal{C}|$ and $\mathcal{P} \subseteq \mathcal{C}'$.
These two conditions on $\mathcal{C}'$ are sufficient because then $\mathcal{C}' \setminus \mathcal{P}$ has size $\coboxicity(K) - |\mathcal{P}|$, and as $E(K') \cap E(C) = \emptyset$ for all $C \in \mathcal{P}$, it follows that $\mathcal{C}' \setminus \mathcal{P}$ is a co-interval cover of $K'$.

\medskip

\noindent
\textbf{Claim 1.} \emph{If $C \in \mathcal{C}(H^-)$, then $E(C) \subseteq E(H^{+})$.}

\smallskip

\noindent
\textit{Proof of Claim.}
Let $u_1u_2 \in E(C) \cap E(H^{-})$, and suppose for a contradiction that there exists $u_3u_4 \in E(C) \setminus E(H^{+})$.
As $v$ separates the two edges into different components, $\{u_1,u_2,u_3,u_4\}$ induces $2K_2$ in $C$, contradicting Lemma~\ref{lem:co-interval-forbidden}.
\hfill $\diamond$

\medskip

\noindent
\textbf{Claim 2.} \emph{$\mathcal{C}(H^{-})$ is a co-interval cover of $H^{-}$, and $|\mathcal{C}(H^{-})| \in \{t^{-}, t^{+}\}$.
Moreover, if $t^{-} = t^{+}$, then $\mathcal{C}$ may be chosen so that $\mathcal{C}(H^{-}) = \cC_{H^+}$.}

\smallskip

\noindent
\textit{Proof of Claim.}
As $t^+ \leq t^- + 1$ by Lemma~\ref{lem:t-sandwich}, we prove $t^- \leq |\mathcal{C}(H^-)| \leq t^+$.
For $|\mathcal{C}(H^-)| \geq t^-$, it is sufficient to show that $\{C[V(H^{-})] : C \in \mathcal{C}(H^{-})\}$ is a co-interval cover of $H^-$.
As co-interval graphs are hereditary, for each $C \in \mathcal{C}(H^{-})$, $C[V(H^{-})]$ is co-interval.
Moreover, since every edge of $H^-$ is covered by some element of $\mathcal{C}(H^-)$, $\{C[V(H^{-})] : C \in \mathcal{C}(H^{-})\}$ covers $E(H^{-})$.
For $|\mathcal{C}(H^-)| \leq t^+$, define $\mathcal{C}^* :=(\mathcal{C} \setminus \mathcal{C}(H^{-})) \cup \cC_{H^+}$.
By Claim~1, every $C \in \mathcal{C}(H^{-})$ has $E(C) \subseteq E(H^{+})$, so $\mathcal{C}^*$ is a co-interval cover of $K$.
Hence, $|\mathcal{C}^*| \le |\mathcal{C}| - |\mathcal{C}(H^{-})| + t^{+}$.
The minimality of $\mathcal{C}$ yields $|\mathcal{C}| \leq |\mathcal{C}^*|$, and so altogether, we get $|\mathcal{C}(H^-)| \leq t^+$.
If $t^{-} = t^{+}$, both bounds coincide, forcing $|\mathcal{C}(H^{-})| = t^{+}$, and so $|\mathcal{C}^*| = |\mathcal{C}|$.
Then, $\mathcal{C}^*$ is itself a minimum-size co-interval cover of $K$ with $\mathcal{C}^*(H^-) = \cC_{H^+}$ by construction, and so we may replace $\mathcal{C}$ by $\mathcal{C}^*$.
\hfill $\diamond$

\medskip

\noindent
\textbf{Claim 3.} \emph{For any $C \in \mathcal{C}$, $C$ is contained in a block ant extending $B$ if and only if $C \in \mathcal{C}(H^{-})$.}

\smallskip

\noindent
\textit{Proof of Claim.}
Suppose $C \in \mathcal{C}(H^-)$.
Then $E(C) \cap E(H^-)$ has an edge $e$.
As $C$ is connected, Lemma~\ref{lem:unique-ant} yields a unique substantial block $Q$ such that $C$ is a subgraph of a block ant extending $Q$.
Suppose for a contradiction that $Q \neq B$.
Note that the blocks in $H$, except for $B$, are not substantial.
As $Q$ is substantial, $V(Q)$ can intersect $V(H)$ only at the anchor $v$.
Thus, $E(C) \cap E(H) \subseteq \delta_K(v)$.
But then $e \in E(C) \cap E(H^-) \subseteq E(C) \cap E(H) \subseteq \delta_K(v)$, and since $e$ is also in $E(H^-)$, it follows that $e \in \delta_K(v) \cap E(H^-) = \emptyset$, a contradiction.

Conversely, suppose $C \notin \mathcal{C}(H^-)$.
Suppose for a contradiction that $C$ is in a block ant extending the unique substantial block $B$.
As $E(C) \cap E(H^-) = \emptyset$ and $E(B) \setminus E(H^-) \subseteq \delta_K(v)$, each edge of $C$ is incident to $v$, so $C$ is a star.
As $K$ has an anchor, it is not a star, so the maximality of $C$ contradicts Proposition~\ref{prop:maximal-star-or-more}.
\hfill $\diamond$

\medskip

Given the claims, we construct $\mathcal{C}'$ for each case.
Throughout, we apply $t^- \leq t \leq t^+$ from Lemma~\ref{lem:t-sandwich}.

\smallskip

\noindent
\textbf{Case 2(a): $t^{-} = t^{+}$.}
By Claim~2, choose $\mathcal{C}$ so $\mathcal{C}(H^{-}) = \cC_{H^+} = \mathcal{P}$.
Set $\mathcal{C}' := \mathcal{C}$.

\smallskip

\noindent
\textbf{Case 2(b): $t^{-} = t < t^{+}$.}
By Claim~2, $|\mathcal{C}(H^-)| \in \{t^-,t^+\} = \{t,t^+\}$.
Suppose $|\mathcal{C}(H^-)| = t^+$.
As $t^{+} \le t^{-}+1 = t+1$ and $t^+ > t$, we get $t^{+} = t+1$.
By Claim~1 and the definition of $\mathcal{C}(H^-)$, every edge covered by $\mathcal{C}(H^-)$ lies in $E(H^+)$.
By the definition of $H^+$, $E(H^+) = E(H) \cup \delta_K(v)$, and so
\[
\mathcal{C}' := (\mathcal{C} \setminus \mathcal{C}(H^-)) \cup \cC_H \cup \{(N_K[v],\delta_K(v))\}
\]
is a co-interval cover of $K$ of size $|\mathcal{C}| - t^{+} + t + 1 = |\mathcal{C}|$, with $\mathcal{P} = \cC_H \subseteq \mathcal{C}'$.
Suppose $|\mathcal{C}(H^-)| = t$. 
In this case, we show that
$$\mathcal{C}' := (\mathcal{C} \setminus \mathcal{C}(H^-)) \cup \mathcal{C}_H$$ 
covers $K$ and satisfies $|\mathcal{C}'| = |\mathcal{C}|$. 
Note that $\mathcal{P} = \mathcal{C}_H \subseteq \mathcal{C}'$.
As $\mathcal{C}'$ replaces $\mathcal{C}(H^-)$ with $\mathcal{C}_H$ and $|\mathcal{C}(H^-)| = |\mathcal{C}_H| = t$, it is sufficient to prove that $\mathcal{C}(H^-)$ covers only edges of $E(H)$.
Suppose for a contradiction that there is an edge $f \notin E(H)$ that is covered by some co-interval graph $C_f \in \mathcal{C}(H^-)$.
Then by Claim~1, $f \in E(H^+)$.
Thus, $f \in E(H^+) \setminus E(H) = \delta_K(v)\setminus E(H)$.
By Claim~3, $C_f$ is contained in a block ant extending $B$.
Moreover, as $E(B)\subseteq E(H)$ and $f \in \delta_K(v)\setminus E(H)$, $v$ must be a hub of this block ant.
As $C_f$ is a maximal connected co-interval subgraph of $K$, by Lemma~\ref{lem:hub-saturation}, $\delta_K(v) \subseteq E(C_f)$.
So, $\mathcal{C}(H^-)$ covers the edge set $E(H)\cup\delta_K(v) = E(H^{+})$.
Thus, $\mathcal{C}(H^-)$ is a co-interval cover of $H^+$ with size $t < t^+$. 
This contradicts the minimality of $t^+$.

\smallskip

\noindent
\textbf{Case~2(c): $t^{-} < t = t^{+}$.}
By Claim~2, $|\mathcal{C}(H^-)| \in \{t^-,t^+\} = \{t^-,t\}$.
Suppose $|\mathcal{C}(H^-)| = t$. 
As $t = t^{+} \le t^{-}+1$ and $t > t^-$, $t = t^{-}+1$.
By Claim~1, every edge covered by $\mathcal{C}(H^-)$ lies in $E(H^+)$.
As $E(H^{+}) = E(H^-)\cup\delta_K(v)$,
\[
\mathcal{C}' := (\mathcal{C}\setminus\mathcal{C}(H^-)) \cup \cC_{H^-} \cup \{(N_K[v],\delta_K(v))\}
\]
is a co-interval cover of $K$ of size $|\mathcal{C}| - t + t^{-}+1 = |\mathcal{C}|$, with $\mathcal{P} = \cC_{H^-} \subseteq \mathcal{C}'$.

Suppose $|\mathcal{C}(H^-)| = t^-$.
Consider the family
$$\mathcal{C}' := (\mathcal{C}\setminus \mathcal{C}(H^-)) \cup \mathcal{C}_{H^-}.$$
As $|\mathcal{C}(H^-)| = |\mathcal{C}_{H^-}|$, $|\mathcal{C}'| = |\mathcal{C}|$. 
Moreover, $\mathcal{P} = \cC_{H^-} \subseteq \mathcal{C}'$.
We now show that $\mathcal{C}'$ covers $K$.
Since $\coboxicity(H) = t = t^- + 1 > t^-$, the set of co-interval graphs $\{C[V(H)]: C \in \mathcal{C}(H^-)\}$ cannot cover all of $E(H)$.
So, some $e = vy_0 \in \delta_H(v)$ (with $y_0 \in V(B) \setminus \{v\})$ is not covered by $\mathcal{C}(H^-)$.
Since $\mathcal{C}$ covers $K$, some $C_e \in \mathcal{C}\setminus \mathcal{C}(H^-)$ must contain $e$.
As $C_e$ is a maximal connected co-interval subgraph of $K$, by Lemma~\ref{lem:unique-ant}, there is a unique substantial block $B'$ such that $C_e$ is contained in a block ant $B'_{u,w}$ for some $u,w \in V(B')$.
By Claim~3, $B' \neq B$.
Further, as $y_0 \notin V(B')$, $v$ is a hub of $B'_{u,w}$, and so $y_0$ is a pendant of $B'_{u,w}$.
Applying Lemma~\ref{lem:hub-saturation}, we obtain $\delta_K(v) \subseteq E(C_e)$.
Since, by Claim~1, every edge covered by $\mathcal{C}(H^-)$ lies in $E(H^+) = E(H^-) \cup \delta_K(v)$, the family $\{C_e\} \cup \mathcal{C}_{H^-} \subseteq \mathcal{C}'$ covers every edge covered by $\mathcal{C}(H^-)$.
Thus, $\mathcal{C}'$ covers $K$.
\end{proof}
\begin{lemma}
\label{lem:iter-optimal}
$|\mathcal{P} \cup \mathrm{Cov}(\Gamma')| \le \coboxicity(\Gamma)$.
\end{lemma}
\begin{proof}
As $B$ is substantial, if $B$ is an edge block, $H^-$ is a star with an edge, and otherwise, $B-v$ contains an edge.
Thus, $E(H^-) \neq \emptyset$, and as $E(H^-) \subseteq E(K)$, $K'$ has strictly fewer edges than $K$. 
So, $\Gamma' = K' \cup (\Gamma\setminus V(K))$ has fewer than $N$ edges.
By the inductive hypothesis, $\mathrm{Cov}(\Gamma')$ is a minimum-size co-interval cover of $\Gamma'$, and so $|\mathrm{Cov}(\Gamma')| = \coboxicity(\Gamma')$.
As $V(K')$ and $V(\Gamma \setminus V(K))$ are disjoint, Lemma~\ref{lem:box_and_disjoint_union} gives $\coboxicity(\Gamma') = \coboxicity(K') + \coboxicity(\Gamma\setminus V(K))$.
So,
\[
|\mathcal{P} \cup \mathrm{Cov}(\Gamma')| = |\mathcal{P}| + \coboxicity(K') + \coboxicity(\Gamma\setminus V(K)).
\]
By Lemma~\ref{lem:iter-decreases}, $\coboxicity(K') \le \coboxicity(K) - |\mathcal{P}|$, so, using $\Gamma = K \cup (\Gamma\setminus V(K))$ and Lemma~\ref{lem:box_and_disjoint_union} once more,
\begin{align}
    |\mathcal{P} \cup \mathrm{Cov}(\Gamma')| \le \coboxicity(K) + \coboxicity(\Gamma\setminus V(K)) = \coboxicity(\Gamma). 
    \tag*{\qed}
\end{align}
\let\qed\relax
\end{proof}
By Lemma~\ref{lem:iter-covers}, $\mathcal{P}\cup\mathrm{Cov}(\Gamma')$ is a co-interval cover of $\Gamma$, so $|\mathcal{P}\cup\mathrm{Cov}(\Gamma')| \ge \coboxicity(\Gamma)$ by Lemma~\ref{lemma-boxicityIntervalOrderCover}.
Combined with Lemma~\ref{lem:iter-optimal}, $|\mathcal{P}\cup\mathrm{Cov}(\Gamma')| = \coboxicity(\Gamma)$, so $\mathcal{P}\cup\mathrm{Cov}(\Gamma')$ is a minimum-size co-interval cover of $\Gamma$, completing the inductive step.
Setting $\Gamma = G$ gives Theorem~\ref{thm:nice-block-algorithm}.

\section{Applications of Theorem~\ref{thm:nice-block-algorithm}}
\label{sec:applications}
In this section, we verify that certain natural families of $2$-connected graphs are tractable, and deduce the corresponding corollaries of Theorem~\ref{thm:nice-block-algorithm}.

\subsection{Proof of Corollary~\ref{cor:tree-coboxicity}}
We begin with the simplest tractable family, the empty family $\mathcal{N} = \emptyset$.
Property (T1) holds vacuously for $\mathcal{N}$.
Since $\mathcal{N} \cup \{K_2\} = \{K_2\}$, verifying (T2) for $\mathcal{N}$ only requires handling pendant augmentations of $K_2$.
Every pendant augmentation of $K_2$ is an induced subgraph of some ant, and since ants are co-interval graphs (recall Figure~\ref{fig:ant-co-interval-representation}), they are also co-interval graphs.
So, a pendant augmentation of $K_2$ forms its own minimum-size co-interval cover.
We record this fact for later use.
\begin{observation}
\label{obs:K2-tractable}
Property (T2) holds when $Q = K_2$.
\end{observation}
Hence, the empty family is tractable.
Let $T$ be a tree.
If $|V(T)| = 1$, then $E(T) = \emptyset$ and $\coboxicity(T) = 0$.
Otherwise, since any block with at least three vertices contains a cycle, every block of $T$ is $K_2$.
Hence every block of $T$ belongs to $\mathcal{N} \cup \{K_2\}$, and Theorem~\ref{thm:nice-block-algorithm} computes $\coboxicity(T)$ in polynomial time.

\subsection{Proof of Corollary~\ref{cor:bounded-block-coboxicity}}
For $k \geq 3$, let $\mathcal{N}_k$ be the family of $2$-connected graphs with order at most $k$ and at least $3$.
Property (T1) is immediate for $\mathcal{N}_k$, so the main task is to verify (T2).
We first show that it suffices to handle pendant augmentations with at most one pendant per vertex, and then solve these by exhaustive search.
Vertices $u$ and $v$ are \emph{twins} in a graph $G$ if $N_G(u) = N_G(v)$.
We state an easy consequence of the definitions of twins and co-interval representations.
\begin{observation}\label{obs:twins-extending-co-interval-representation}
    If vertices $u$ and $v$ are twins in $H$, and $I$ is a co-interval representation of $H-u$, then setting $I(u) := I(v)$, the family $I \cup \{I(u)\}$ is a co-interval representation of $H$.
\end{observation}
\begin{lemma}
\label{lem:extra-pendant-invariant}
Let $G$ be a graph containing a vertex $v$ with pendant $\ell$, and $\ell'$ a vertex outside $V(G)$.
Let $\mathcal{C} = \{H_1, H_2, \ldots, H_t\}$ be a co-interval cover of $G$ such that $\ell v \in E(H_1)$.
Then $\{H_1 \cup (\{\ell',v\},\{\ell'v\}),$ $H_2,$ $\ldots,$ $H_t\}$ is a co-interval cover of $G \cup (\{\ell',v\},\{\ell'v\})$ obtained from $\mathcal{C}$ in constant time.
\end{lemma}
\begin{proof}
    Let $I_1$ be a co-interval representation of $H_1$. 
    Set $H_1^+ := H_1 \cup (\{\ell',v\},\{\ell'v\})$ and $G^+ := G \cup (\{\ell',v\},\{\ell'v\})$.
    As $\ell$ and $\ell'$ are pendants of $v$ in $G^+$ and $H_1^+$ is a subgraph of $G^+$ containing $\ell v$ and $\ell' v$, $\ell$ and $\ell'$ are twins in $H_1^+$.
    By Observation~\ref{obs:twins-extending-co-interval-representation}, a co-interval representation of $H_1^+$ is obtained from $I_1$ by adding $I_1(\ell') := I_1(\ell)$.
    Thus $\{H_1^+, H_2, \ldots, H_t\}$ is a co-interval cover of $G^+$, which is obtained from $\mathcal{C}$ by adding one vertex and one edge to $H_1$, which takes constant time.
\end{proof}

\begin{lemma}
\label{lem:reduce-to-single-pendant}
Let $Q$ be a graph, and let $q$ be a non-decreasing function such that, for every $x' : V(Q) \to \{0,1\}$, a minimum-size co-interval cover of $Q^{x'}$ can be computed in time $q(|E(Q^{x'})|)$.
Then, for every $x : V(Q) \to \mathbb{Z}_{\geq 0}$, a minimum-size co-interval cover of $Q^{x}$ can be computed in time $O\big(q(|E(Q^{x'})|) + |E(Q^{x})|\big)$.
\end{lemma}
\begin{proof}
Let $x'_v = \min(x_v, 1)$ for each $v \in V(Q)$.
Since $Q^{x'}$ is a subgraph of $Q^{x}$, we can compute a minimum-size co-interval cover $\mathcal{C}'$ of $Q^{x'}$ in time $q(|E(Q^{x'})|)$.
The graph $Q^{x}$ is obtained from $Q^{x'}$ by attaching $x_v - 1$ further pendants to each $v$ with $x_v \geq 2$.
By repeated application of Lemma~\ref{lem:extra-pendant-invariant}, all of them can be added to a single piece of $\mathcal{C}'$ containing the edge between $v$ and its pendant in $Q^{x'}$, and the result is a minimum-size co-interval cover of $Q^{x}$.
These pieces are found by a single scan of $\mathcal{C}'$, whose cost is at most the time $q(|E(Q^{x'})|)$ needed to output $\mathcal{C}'$, and each of the at most $|E(Q^{x})|$ additions takes constant time.
\end{proof}
\begin{lemma}
\label{lem:bounded-order-tractable}
For every $k \geq 3$, the family $\mathcal{N}_k$ is tractable.
Moreover, for every $Q \in \mathcal{N}_k \cup \{K_2\}$ and $x : V(Q) \to \mathbb{Z}_{\geq 0}$, a minimum-size co-interval cover of $Q^{x}$ can be computed in time $O\big(k^{k^2 + 6k + 7} + |E(Q^x)|\big)$.
\end{lemma}
\begin{proof}
For (T1), fix $Q \in \mathcal{N}_k$ and $u \in V(Q)$.
Since $Q$ is $2$-connected with at least three vertices, $Q - u$ is connected and has at least two vertices, so each of its blocks is $2$-connected and has fewer than $k$ vertices.
Hence each block of $Q - u$ is either $K_2$ or belongs to $\mathcal{N}_k$.
For (T2), fix $Q \in \mathcal{N}_k \cup \{K_2\}$.
The case $Q = K_2$ follows by Observation~\ref{obs:K2-tractable}, 
so we assume $Q \in \mathcal{N}_k$.
Let $x: V(Q) \to \{0,1\}$.
So, by Lemma~\ref{lem:reduce-to-single-pendant}, it suffices to show that a minimum-size co-interval cover of $Q^{x}$ can be computed in time $O\big(k^{k^2 + 6k + 7}\big)$.
As $Q$ is a block of order at least $3$, it is not a star, and so by Lemma~\ref{lem:unique-ant}, the edge set of any maximal co-interval subgraph of $Q^x$ is a subgraph of a unique block ant that extends a substantial block. 
As $Q$ is the only substantial block of $Q^x$, such a block ant is of the form $Q_{u,v}$ for some $u,v \in V(Q)$, where $u$ and $v$ each have at most one neighbour outside $V(Q)$.
By~\cite[Corollary~1]{2023_Caoduro_GD}, every graph $G$ has at most $|V(G)|!$ maximal co-interval subgraphs, each constructible in time $O(|V(G)|^{2})$.
So, since $Q_{u,v}$ has order at most $k+2$, for each pair $u,v \in V(Q)$, there are at most $(k+2)!$ maximal co-interval subgraphs of $Q_{u,v}$.
Thus, ranging over the at most $k^2$ block ants in $Q^x$, the family $\mathcal{M}$ of maximal co-interval subgraphs of $Q^x$ satisfies $|\mathcal{M}| \leq k^2 (k+2)!$, which can be computed in time $O(k^2(k+2)! \cdot (k+2)^2)$.

Every graph on $n$ vertices has co-boxicity at most $\lfloor n/2 \rfloor$~\cite[Corollary~7.1]{Rob}.
This implies that $\coboxicity(Q^{x}) \leq k$, and so a minimum-size co-interval cover of $Q^x$ has at most $k$ elements.
Hence, to obtain a minimum-size cover consisting of elements of $\mathcal{M}$, it suffices to check for each of the at most
$$\binom{k^2(k+2)!}{\leq k} \leq \big(k^2(k+2)!\big)^{k} \leq k^{2k}(k+2)^{k(k+2)}$$
subsets of $\mathcal{M}$ of size at most $k$, whether it covers $E(Q^{x})$, and to return the smallest covering subset.
Since $|V(Q)| \leq k$, we have $|V(Q^{x})| \leq 2k$ and $|E(Q^{x})| \leq \binom{k}{2} + k \leq k^2$.
So, each check takes time $O(k \cdot |E(Q^{x})|) = O(k^3)$.
Thus, the search takes time $O\big(k^3 \cdot k^{2k}(k+2)^{k(k+2)}\big)$, which dominates the time needed to compute $\mathcal{M}$.
As $k+2 \leq k^{(k+2)/k}$ for $k \geq 3$, $(k+2)^{k(k+2)} \leq k^{(k+2)^2} = k^{k^2 + 4k + 4}$.
Hence the running time is $O\big(k^{k^2 + 6k + 7}\big)$.
\end{proof}
Let $G$ be a graph where each block has at most $k$ vertices.
Then every block of $G$ is in $\mathcal{N}_k \cup \{K_2\}$.
By Lemma~\ref{lem:bounded-order-tractable}, Theorem~\ref{thm:nice-block-algorithm} applies with $\mathcal{N} = \mathcal{N}_k$ and $q(|E(G)|) = O\big(k^{k^2 + 6k + 7} + |E(G)|\big)$, and so $\coboxicity(G)$ is computable in time
$$
O\big(|E(G)|^2 \cdot (k^{k^2 + 6k + 7} + |E(G)|)\big) = O\big(k^{k^2 + 6k + 7} \cdot |E(G)|^2 + |E(G)|^3\big),
$$
which proves Corollary~\ref{cor:bounded-block-coboxicity}.

\subsection{Proof of Corollary~\ref{cor:block-graph-coboxicity}
}
Lemma~\ref{lem:ECS_are_subgraph_of_big_ant} shows that every connected co-interval subgraph fits inside some block ant, and while the containing block ant is not co-interval in general, if the block is a clique, we will show that it is a co-interval graph.
A \emph{clique ant} is a block ant where the block is a clique.
\begin{example}\label{eg:ECS_example}
    Figure~\ref{fig:ECS_example} shows a graph $G$ with clique ant subgraph $Q_{u,v}$ and corresponding co-interval representation $I$.
    Notice that all intervals corresponding to vertices in $Q$ are pairwise disjoint, and an interval corresponding to any pendant in $Q_{u,v}$ intersects all other intervals except the one corresponding to the pendant's neighbour.
    \begin{figure}
        \centering
        \includegraphics[scale=0.85]{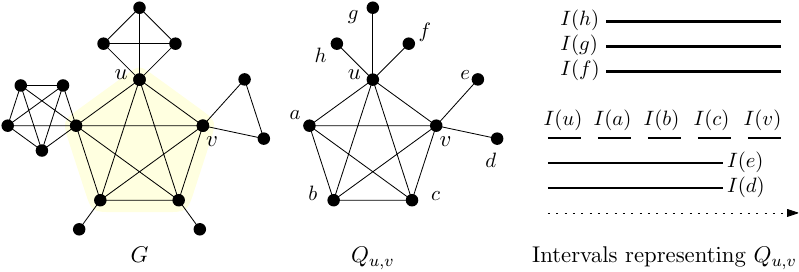}
        \caption{
        A clique ant subgraph with corresponding co-interval representation.
        }
        \label{fig:ECS_example}
    \end{figure}
\end{example}
The construction in the proof below generalizes the co-interval representation used in Example~\ref{eg:ECS_example} to every clique ant.
\begin{lemma}\label{lem:ECS_are_interval_order}
    Let $Q$ be a clique block in a graph $G$, and let $u,v \in V(Q)$.
    Then $Q_{u,v}$ is a co-interval graph.
\end{lemma}
\begin{proof}
    Let $\mathcal{I} = \{I(t) : t \in V(Q)\}$ be a family of pairwise disjoint closed intervals where $I(u)$ is the leftmost interval and, if $u \neq v$, $I(v)$ is the rightmost interval in the family.
    For each vertex $w \in N_G(u) \setminus V(Q)$ (respectively, $w' \in N_G(v) \setminus V(Q)$), add $I(w)$ (respectively, $I(w')$) to $\mathcal{I}$ intersecting all intervals in $\mathcal{I}$ except $I(u)$ (respectively, $I(v)$).
    For $w,w' \in V(Q_{u,v})$, $I(w) \cap I(w') = \emptyset$ if and only if $ww' \in E(Q_{u,v})$.
    Thus, $I : V(Q_{u,v}) \xrightarrow{} \mathcal{I}$ is a co-interval representation of $Q_{u,v}$.
    %
    %
\end{proof}
\begin{lemma}
\label{lem:whiskered-clique-coboxicity}
Let $Q$ be a clique on at least two vertices, let $x : V(Q) \to \mathbb{Z}_{\geq 0}$, and let $r$ be the number of vertices $v \in V(Q)$ with $x_v \geq 1$.
Then $\coboxicity(Q^{x}) = \max(1, \lceil r/2 \rceil)$, and a minimum-size co-interval cover of $Q^{x}$ consisting of clique ants can be computed in time polynomial in $|E(Q^{x})|$.
\end{lemma}
\begin{proof}
    As $E(Q) \neq \emptyset$, $\coboxicity(Q^x) \geq 1$.
    If $r \leq 2$, then $Q^x = Q_{u,v}$ for some $u,v \in V(Q)$, and so by Lemma~\ref{lem:ECS_are_interval_order}, $\{Q\}$ is a co-interval cover of $Q^x$, implying $\coboxicity(Q^x) = 1$.
    Suppose $r \geq 3$.
    We show that $\coboxicity(Q^x)$ is bounded above and below by $\ceil{r/2}$.
    Let $v_1, v_2, \ldots, v_r$ be the vertices in $V(Q)$ with pendants outside $V(Q)$. 
    By Lemma~\ref{lem:ECS_are_interval_order}, $Q_{v_i,v_j}$ is a co-interval subgraph of $Q^x$ for all $1 \leq i \leq j \leq r$.
    As all edges of $Q^x$ are incident to vertices in $V(Q)$, $\{Q_{v_1,v_r}, Q_{v_2,v_{r-1}},$ $\ldots,$ $Q_{v_{\ceil{r/2}},v_{\floor{r/2}+1}}\}$ is a co-interval cover of $Q^x$ of size $\ceil{r/2}$ and it is easily computable in time polynomial in $|E(Q^x)|$.
    Now we show $\coboxicity(Q^x) \geq \ceil{r/2}$.
    Let $\mathcal{C}$ be a minimum-size co-interval cover of $Q^x$, where each element is maximal.
    As $r \geq 3$, $Q^x$ is not a star, and so by Lemma~\ref{lem:unique-ant}, every maximal co-interval subgraph of $Q^x$ is in a unique clique ant $Q_{u,v}$ for some $u,v \in V(Q)$, where $u$ and $v$ have pendants outside $V(Q)$.
    As $Q_{u,v}$ is co-interval, the maximality of the elements in $\mathcal{C}$ implies that $\mathcal{C}$ consists entirely of clique ants extending $Q$.
    The only edges outside $E(Q)$ that each $Q_{u,v}$ in $\mathcal{C}$ can cover are incident to $u$ and $v$, and all such edges are covered by $\mathcal{C}$. 
    So, $\coboxicity(Q^x) = |\mathcal{C}| \geq \ceil{r/2}$.
\end{proof}

\begin{lemma}
\label{lem:cliques-tractable}
The family $\mathcal{N} = \{K_k : k \geq 3\}$ is tractable.
\end{lemma}
\begin{proof}
For (T1), let $Q = K_k$ with $k \geq 3$, and let $u \in V(Q)$.
Then $Q - u = K_{k-1}$ is $2$-connected, and hence a single block, which belongs to $\mathcal{N}$ if $k \geq 4$ and equals $K_2$ if $k = 3$.
Property (T2) follows from Lemma~\ref{lem:whiskered-clique-coboxicity}, which applies to every $Q \in \mathcal{N} \cup \{K_2\} = \{K_k : k \geq 2\}$ and every $x : V(Q) \to \mathbb{Z}_{\geq 0}$.
\end{proof}
Since every block of a block graph belongs to $\{K_k : k \geq 3\} \cup \{K_2\}$, Theorem~\ref{thm:nice-block-algorithm} and Lemma~\ref{lem:cliques-tractable} yield Corollary~\ref{cor:block-graph-coboxicity}.
\subsection{Proof of Corollary~\ref{cor:cactus-coboxicity}}
Let $C_k$ denote the cycle of length $k$.
In this section, we show that the family $\{C_k : k \geq 3\}$ is tractable, which, by Theorem~\ref{thm:nice-block-algorithm}, gives Corollary~\ref{cor:cactus-coboxicity}.

Pendant augmentations of short cycles are already handled by Lemma~\ref{lem:bounded-order-tractable}, since $C_3, C_4 \in \mathcal{N}_4$.
For longer cycles, we show that some minimum-size cover uses the ant of a cycle edge at a vertex with a pendant, and that the rest of the graph is a tree.
We first prove a general proposition about graphs with no induced $C_3$ and $C_4$.
\begin{proposition}\label{prop:girth-at-least-5-co-intervals}
    Let $G$ be a graph with at least one edge and no cycles of length $3$ or $4$, and let $H$ be a maximal connected co-interval subgraph of $G$.
    Then $H$ is the $(u,v)$-ant of some edge $uv$ of $G$.
\end{proposition}
\begin{proof}
    If $H$ is not a tree, then it contains an induced cycle with $t$ vertices.
    Since every cycle with at least $6$ vertices contains an induced $2K_2$ and, by Lemma~\ref{lem:co-interval-forbidden}, $H$ contains no induced $2K_2$ or $C_5$, we have $t < 5$.
    Then as $H$ is a subgraph of $G$, which contains no $C_3$ or $C_4$, $t < 3$. 
    So, $H$ is a tree.
    Moreover, as any path with more than $4$ vertices in a tree also contains an induced $2K_2$, the number of vertices in a longest path in $H$ is at most $4$.
    It follows that $H$ is a subgraph of some $(u,v)$-ant in $G$.
    Since ants are co-interval and $H$ is maximal, $H$ equals this ant.
\end{proof}
\begin{lemma}
\label{lem:whiskered-cycle-reduction}
Let $C$ be a cycle with vertices $v_1, \dots, v_k$ in cyclic order, where $k \geq 5$, and let $x : V(C) \to \mathbb{Z}_{\geq 0}$ with $x_{v_1} \geq 1$.
For an edge $uv$ of $C$, let $T_{uv} = (C - \{u,v\})^{x}$, with $x$ restricted to $V(C) \setminus \{u,v\}$.
Then
\[
\coboxicity(C^{x}) = 1 + \min\big(\coboxicity(T_{v_kv_1}),\ \coboxicity(T_{v_1v_2})\big),
\]
and a minimum-size co-interval cover of $C^{x}$ can be computed in time polynomial in $|E(C^{x})|$.
\end{lemma}
\begin{proof}
For an edge $uv$ of $C$, let $A_{uv}$ be the $(u,v)$-ant of $C^{x}$, which is co-interval.
Its edge set is $\delta(u) \cup \delta(v)$, so $E(C^{x})$ is the disjoint union of $E(A_{uv})$ and $E(T_{uv})$.
For the upper bound, as $E(C^x) = E(A_{uv}) \cup E(T_{uv})$, adding $A_{uv}$ to a minimum-size co-interval cover of $T_{uv}$ gives a co-interval cover of $C^{x}$, so $\coboxicity(C^{x}) \leq 1 + \coboxicity(T_{uv})$.
Each $T_{uv}$ is a tree, so Algorithm~\ref{alg:coboxicity} with $\mathcal{N} = \emptyset$ computes a minimum-size co-interval cover of it in polynomial time.
Doing this for $uv \in \{v_kv_1, v_1v_2\}$ and keeping the smaller of the two resulting covers gives a cover of the claimed size.
For the lower bound, let $\mathcal{C}$ be a minimum-size co-interval cover of $C^{x}$, suppose each element is maximal, and let $H \in \mathcal{C}$ contain a pendant edge $v_1\ell$.
As $k \geq 5$, $C^x$ contains no $C_3$ or $C_4$, and so by Proposition~\ref{prop:girth-at-least-5-co-intervals}, $H = A_{u,v}$ for some $u,v \in V(C^x)$.
Since $H$ contains $v_1\ell$, we have that $H = A_{v_1w}$ with $w \in \{v_2,v_k\}$.
Since $E(C^x) = E(A_{v_1w}) \cup E(T_{v_1w})$, and $E(A_{v_1w}) \cap E(T_{v_1w}) = \emptyset$, the elements of $\mathcal{C}\setminus \{A_{v_1w}\}$ cover $E(T_{v_1w})$, and so $\coboxicity(T_{v_1w}) \leq |\mathcal{C}| - 1 = \coboxicity(C^x) - 1$, giving the lower bound.
\end{proof}
\begin{lemma}
\label{lem:cycles-tractable}
The family $\mathcal{N} = \{C_k : k \geq 3\}$ is tractable.
\end{lemma}
\begin{proof}
For (T1), let $Q = C_k$ with $k \geq 3$, and let $u \in V(Q)$.
Then $Q - u$ is a path, all of whose blocks are $K_2$.
For (T2), let $Q \in \mathcal{N} \cup \{K_2\}$ and $x : V(Q) \to \mathbb{Z}_{\geq 0}$.
If $Q = K_2$, we use Observation~\ref{obs:K2-tractable}, and if $Q \in \{C_3, C_4\} \subseteq \mathcal{N}_4$, we use Lemma~\ref{lem:bounded-order-tractable}.
So let $Q = C_k$ with $k \geq 5$, with vertices $v_1, \dots, v_k$ in cyclic order.
If $x \equiv 0$, then $Q^{x} = C_k$, and the $(v_i, v_{i+1})$-ants for $i \in \{2, 5, 8, \dots\}$, indices modulo $k$, form a co-interval cover of size $\lceil k/3 \rceil$, which is minimum by~\cite{1983_Cozzens}.
Otherwise, we may assume that $x_{v_1} \geq 1$, and Lemma~\ref{lem:whiskered-cycle-reduction} applies.
\end{proof}
Since every block of a cactus graph belongs to $\{C_k : k \geq 3\} \cup \{K_2\}$, Theorem~\ref{thm:nice-block-algorithm} and Lemma~\ref{lem:cycles-tractable} yield Corollary~\ref{cor:cactus-coboxicity}.

\section{Conclusion}
\label{section-discussion}
We proved that for a graph $G$ with all blocks in a tractable family, computing $\coboxicity(G)$ is in \textbf{P}.
To do this, we showed that $\mathrm{Cov}$ successfully processes $G$ recursively to build a minimum-size co-interval cover.
Then, by proving tractability, we exemplified this new method on trees, graphs with bounded block order, block graphs, and cactus graphs.
This motivates future work toward proving tractability for other 2-connected graph families.

A natural next step in developing the method is to determine whether Property (T2) of Definition~\ref{def:tractable} can be relaxed to require a polynomial-time algorithm only for the co-boxicity of blocks, rather than their pendant augmentations.

Co-boxicity is an instance of a broader family of covering parameters studied by Cozzens and Roberts~\cite{1989_Cozzens}.
Our method relies on specific properties of co-interval graphs, so extending it to other such families may require new ideas, though it seems promising whenever the covering family is $2K_2$-free.
One example is the threshold co-dimension, whose covering family is the class of threshold graphs, i.e.\ graphs with no induced $2K_2$, $P_4$, or $C_4$.
A similar technique already yields a polynomial-time algorithm for the threshold co-dimension of block graphs in an earlier version of this work~\cite{2025_Caoduro}.

\section*{Acknowledgements}
The authors are grateful to Amna Adnan, Matthew Barclay, Josh Childs, and the coordinators of the Pacific Institute for the Mathematical Sciences (PIMS) VXML program for initiating this research collaboration.
M. Caoduro was supported by a Natural Sciences and Engineering Research Council of Canada Discovery Grant [RGPIN-2021-02475].
W. Evans was supported by NSERC Discovery Grant RGPIN-2022-04449.
We gratefully acknowledge that this research was supported in part by PIMS.
%
\subsubsection*{Declaration of generative AI in the writing process.}
During the preparation of this work, the authors used Claude Sonnet~5.5 and Claude Opus~5.5 (Anthropic) to proofread the manuscript and to draft some paragraphs in the introduction and in Section~\ref{sec:applications} from prompts provided by the authors.
After using these tools, the authors reviewed and edited the content as needed and take full responsibility for the content of the publication.
\bibliographystyle{splncs04}
\bibliography{references_short}
\end{document}